\documentclass[12pt, a4paper]{amsart}

\usepackage[unicode]{hyperref}
\hypersetup{colorlinks,bookmarksopen,bookmarksnumbered,citecolor=blue,pdfstartview=FitH}

\usepackage{color}

\usepackage[english]{babel}
\usepackage{amssymb,amsmath,amsthm}
\usepackage{verbatim}

\usepackage[all]{xy}

\def\AA{\mathcal A}
\def\BB{\mathcal B}
\def\CC{\mathcal C}
\def\DD{\mathcal D}
\def\EE{\mathcal E}

\def\eps{\varepsilon}
\def\={\;=\;}
\def\bal{\begin{aligned}}
\def\eal{\end{aligned}}
\def\be{\begin{equation}\label}
\def\ee{\end{equation}}

\newcommand{\Vect}{{\rm Vect_k}}
\newcommand{\Hom}{{\rm Hom}}

\newcommand{\id}{{\rm id}}

\title[Group actions on categories and Elagin's Theorem revisited]{Group actions on categories and\\Elagin's Theorem Revisited}

\author{Evgeny Shinder}

\begin{document}

\maketitle

\begin{abstract}
After recalling basic definitions and constructions for a finite group $G$ action on a $k$-linear category we give a concise proof
of the following theorem of Elagin: if $\CC = \langle \AA, \BB \rangle$ is a semiorthogonal decomposition of a triangulated category
which is preserved by the action of $G$, and $\CC^G$ is triangulated, then there is a semiorthogonal decomposition $\CC^G = \langle \AA^G, \BB^G \rangle$.
We also prove that any $G$-action on $\CC$ is weakly equivalent to a strict $G$-action which is the analog of the Coherence Theorem for 
monoidal categories.

\medskip
Keywords: group actions on categories, derived categories of coherent sheaves, Elagin's Theorem

\medskip
Mathematics Subject Classiﬁcation 2010: 
14L30, 
18E30. 

\end{abstract}


\section{Introduction}

\subsection{} The setting of finite groups acting on categories is a well-studied ground, see e.g. \cite{Deligne, Sosna, GK, E1, E2} and
references therein. A useful way to define the action is to require for every $g \in G$ an autoequivalence $\rho_g\colon \CC \to \CC$
together with a choice of isomorphisms $\rho_g \rho_h \simeq \rho_{gh}$ satisfying a cocycle condition, see \ref{def-action}.
One would then study the category of equivariant objects $\CC^G$, see \ref{def-equiv}.

\subsection{}\label{ex-Db} For instance, if $\CC = \DD^b(X)$ is the derived category of coherent sheaves on a variety $X$ then a $G$-action on 
$X$ induces a $G$-action on $\CC$, and furthermore $\CC^G$ can be interpreted as the derived category of coherent sheaves on the quotient stack $X/G$.

\subsection{} The main goal of this paper is to give a direct proof of the Theorem of Elagin \cite{E1,E2} stating that if
$\CC = \langle \AA, \BB \rangle$ is a semi-orthogonal decomposition of triangulated categories and $G$ is
finite group acting on $\CC$ by triangulated autoequivalences in such a way that the category of equivariant objects $\CC^G$ is triangulated
and preserving $\AA$ and $\BB$, then there is a semi-orthogonal decomposition $\CC^G = \langle \AA^G, \BB^G \rangle$,
see Theorem \ref{thm-elagin}.
In the setup of \ref{ex-Db} this Theorem is often quite useful in constructing semiorthogonal decompositions for 
the quotient stack $\DD^b(X/G)$ from semiorthogonal decompositions of $\DD^b(X)$.

\subsection{} In our proof we construct the functors $\CC^G \to \AA^G$ and
$\CC^G \to \AA^G$ adjoint to the inclusion functors.
The key step in the proof is to show that if $\Phi\colon \AA \to \CC$ is a $G$-equivariant functor which admits
a left or right adjoint functor $\Psi$, then $\Psi$ is automatically equivariant: see Proposition \ref{adj-mates}.

\subsection{} We also prove that every $G$-action on a category $\CC$ is $G$-weakly equivalent to a strict $G$-action, that is to
an action satisfying $\rho_g \rho_h = \rho_{gh}$, see Theorem \ref{strict}.
This is analogous to the Coherence Theorem for monoidal categories: every monoidal category is equivalent
to a strict monoidal category, see e.g. \cite[1.2.15]{L}.

\subsection{} In order to formulate and prove these facts we need to develop the language of $G$-functors,
$G$-natural transformations and so on. 
Perhaps relevant definitions and constructions are well-known to experts but we include these
for completeness as we could not find the reference that fits our purpose.

\subsection{} All categories, functors etc are $k$-linear where $\mathrm{char}(k) = 0$. Groups acting on categories are finite
and we denote by $1 \in G$ the neutral element of the group.

We use the symbol ``$\circ$'' to denote vertical composition of natural transformations of functors,
the other types of compositions are denoted by concatenation.

\subsection{Acknowledgements:} We thank A.\,Elagin, S.\,Galkin, N.\,Gurski, T.\,Leinster and F.\,Petit for useful conversations
and e-mail communications and the referee for the suggestions on improving the exposition. 

\section{$G$-categories and equivariant objects}

\subsection{}

\label{def-action}

By a {\bf $G$-action on $\CC$} we mean the following data \cite[Def.\,3.1]{E2}:

\begin{itemize}
\item For each element $g \in G$ an autoequivalence $\rho_g\colon \CC \to \CC$

\item For each pair $g,h \in G$ an isomorphism of functors 
\[
\phi_{g,h}\colon \rho_g \rho_h \cong \rho_{g h}. 
\]


\end{itemize}
The data must satisfy the following associativity axiom: for all $g,h,k \in G$ the diagram of functors $\CC \to \CC$ is commutative:
\[\xymatrix{
\rho_g \rho_h \rho_k \ar[d]_{\phi_{g,h}\rho_k} \ar[r]_{\rho_g \phi_{h,k}} & \rho_g \rho_{h k} \ar[d]_{\phi_{g,hk}} \\
\rho_{gh} \rho_k \ar[r]_{\phi_{gh,k}} &  \rho_{ghk}
}\]



\subsection{}

It follows from the definition that there is an isomorphism of functors
\[
\phi_1\colon \rho_1 \simeq id
\]
obtained by post-composing 
\[
\phi_{1,1}\colon \rho_1 \rho_1 \to \rho_1 
\]
with $\rho_1^{-1}$. That is we have
\[
\phi_{1,1} = \rho_{1} \phi_{1}. 
\]
Furthermore one can show that $\phi_1$ satisfies \cite[2.1.1(e)]{GK}:
\[\bal
\phi_{g,1} &= \rho_g \phi_1\colon \rho_g \rho_1 \to \rho_g \\
\phi_{1,g} &= \phi_1 \rho_g\colon \rho_1 \rho_g \to \rho_g \\
& \\
\eal\]
so that
definition \ref{def-action} coincides with that of \cite[2.1]{GK}.

On the other hand if one asks for $\phi_1$ to be the identity transformation,
one gets a slightly stronger definition of a $G$-descent datum of \cite[Def.\,1.1]{N}.

\subsection{}

Using the language of monoidal functors \cite[Def.\,1.2.10]{L} one can give a very
concise definition of a group acting on a category. For that consider
$G$ as a monoidal category: $G$ is discrete as a category and its monoidal structure defined by
\[
g \otimes h = gh
\]
\[
id_g \otimes id_h = id_{gh}. 
\]

Now a $G$-action on $\CC$ amounts to the same thing as an action of monoidal category $G$ on $\CC$ \cite[Ex.\,1.2.12]{L},
i.e. a weak monoidal functor 
\[
\rho\colon G \to [\CC, \CC] 
\]
where on the right is the category of functors $\CC \to \CC$ with monoidal structure
given by composing functors.

\subsection{}\label{def-equiv}

One defines the {\bf category of $G$ equivariant objects $\CC^G$} \cite{E2,GK} as follows:
objects of $\CC^G$ are linearized objects, i.e. objects $c \in \CC$ equipped with
isomorphisms
\[
\theta_g\colon c \to \rho_g (c), \; g \in G
\]
satisfying the condition that the diagrams are commutative:
\[\xymatrix{
c \ar[d]_{\theta_{gh}} \ar[r]_{\theta_g} & \rho_g(c) \ar[d]_{\rho_g \theta_h} \\
\rho_{gh}(c) & \ar[l]_{\phi_{g,h}(c)} \rho_g (\rho_h (c))
}\]

Morphisms of equivariant objects consist of those morphisms of the
underlying objects in $\CC$ 
which commute with all $\theta_g$, $g \in G$.

\section{$G$-functors and $G$-natural transformations}

\subsection{}\label{pentagon}

Given two categories $\CC$, $\DD$ with $G$-actions and a functor $\Phi\colon \CC \to \DD$, $\Phi$ is called 
a {\bf right lax $G$-functor} if
there are given natural transformations 
\[
\delta_g\colon \rho_g \Phi \to \Phi \rho_g
\]
such that the two natural transformations $\rho_g \rho_h \Phi \to \Phi \rho_{gh}$ coincide:
\[\xymatrix{
 & & \rho_{gh} \Phi \ar[drr]^{\delta_{gh}} & &  \\
\rho_g\rho_h\Phi \ar[urr]^{\phi_{g,h} \Phi} \ar[dr]_{\rho_g \delta_h} &  &  &  & \Phi\rho_{gh} \\
& \rho_g\Phi\rho_h \ar[rr]_{\delta_g \rho_h} &  &  \Phi\rho_{g}\rho_h \ar[ur]_{\Phi \phi_{g,h}} & & \\
}\]
This commutative diagram is called the pentagon axiom.

\medskip


Similarly $\Phi$ is called a {\bf left lax $G$-functor} if there are given
natural transformations 
\[
\delta_g\colon \Phi \rho_g \to \rho_g \Phi 
\]
satisfying the dual pentagon axiom.

A right (or left) lax $G$-functor $\Phi$ is called a {\bf weak $G$-functor} 
if all $\delta_g$ are isomorphisms.

\medskip

The following lemma is a useful criterion for a weak $G$-functor.

\subsection{Lemma}
\label{lem-weak}
Let $\Phi$ be a right (or left) lax $G$-functor.
The following conditions are equivalent:

\begin{enumerate}
 \item The natural transformation $\delta_1\colon \rho_1 \Phi \to \Phi \rho_1$ is an isomorphism.
 \item $\Phi$ satisfies the identity element axiom:
\[
\Phi \phi_1 \circ \delta_1 = \phi_1 \Phi\colon \rho_1 \Phi \to \Phi. 
\]
 \item $\Phi$ is a weak $G$-functor.
\end{enumerate}

\subsection{Proof}

Implications $(3) \implies (1)$, $(2) \implies (1)$ are obvious.
Let us prove that $(1) \implies (3)$.
Consider the case of the right lax $G$-functor. Applying the pentagon axiom 
to the pair $(g^{-1},g)$ gives:
\[
\delta_{g^{-1}} \rho_g \circ \rho_{g^{-1}} \delta_g = \Phi \phi_{g^{-1},g}^{-1} \circ \delta_1 \circ \phi_{g^{-1},g} \Phi. 
\]
Since the natural transformation on the right-hand side is an isomorphism (note that $\delta_1$ is an isomorphism
by the identity element axiom) and $\rho_{g}$, $\rho_{g^{-1}}$ are
equivalences, it follows that
$\delta_{g^{-1}}$ is left invertible and $\delta_g$ is right invertible.
Thus we see that all $\delta_g$ are isomorphisms.

\medskip

Now we prove $(1) \implies (2)$. Consider the natural transformation
\[
\eps = \Phi \phi_1 \circ \delta_1 \circ \phi_1^{-1} \Phi\colon \Phi \to \Phi.
\]
We are given that $\eps$ is an isomorphism and we need to prove that $\eps$ is in fact an identity.


We use Lemma \ref{G-compos} applied to the trivial group $H := \{1\}$ and the composition
\[
(\CC, id) \underset{(id, \phi_1)}{\to} (\CC, \rho_1) \underset{(\Phi, \delta_1)}\to (\DD, \rho_1) \underset{(id, \phi_1^{-1})}\to (\DD, id) 
\]
which gives a lax $G$-functor
\[
(\CC, id) \underset{(\Phi, \eps)}{\to} (\DD, id).
\]
The pentagon axiom for this functor yields
\[
\eps^2 = \eps 
\]
and we deduce that $\eps = id$.

\subsection{Lemma}
\label{G-compos}

If $(\Phi, \delta^\Phi)\colon \CC \to \DD$, $(\Psi, \delta^\Psi)\colon \DD \to \EE$ are right/left/weak $G$-functors, then their composition
$(\Psi \Phi, \Phi \delta^{\Psi} \circ \delta^\Phi \Psi)$ is a right/left/weak $G$-functor.

\medskip

For the proof one needs to check that the composition satisfies the pentagon and/or the identity element axioms; this is a straightforward check.

\subsection{Lemma}

A weak $G$-functor $\Phi\colon \CC \to \DD$ induces a functor on the categories of equivariant objects 
\[
\Phi^G\colon \CC^G \to \DD^G
\]
such that the following diagram is commutative
\[\xymatrix{
\CC^G \ar[d] \ar[r]^{\Phi^G} & \DD^G \ar[d] \\
\CC \ar[r]^{\Phi} &  \DD
}\]

\subsection{Proof}
For $(c, \theta) \in \CC^G$ we define linearization on $\Phi(c)$ as a composition of isomorphisms
\[
\Phi(c) \to \Phi \rho_g(c) \to \rho_g \Phi(c)
\]
of $\Phi \theta_g$ with $\delta_g$. It is now a standard check that $\Phi(c)$ becomes an equivariant object
and that $\Phi^G$ is a functor.

\subsection{Definition} A natural transformation between two weak $G$-functors $\mu\colon \Phi_1 \to \Phi_2\colon \CC \to \DD$
is called a $G$-natural transformation if for every $g \in G$ the following diagram commutes:
\[\xymatrix{
\rho_g \Phi_1 \ar[d]_{\delta_{1,g}} \ar[r]^{\rho_g \mu} & \rho_g \Phi_2 \ar[d]_{\delta_{2,g}} \\
\Phi_1 \rho_g \ar[r]^{\mu \rho_g} &  \Phi_2 \rho_g.
}\]

\subsection{Lemma}
\label{lem-Gnat}
A $G$-natural transformation $\mu$ between two weak $G$-functors $\Phi_1, \Phi_2\colon \CC \to \DD$
induces a natural transformation $\mu^G\colon \Phi_1^G \to \Phi_2^G$.

\subsection{Proof}

To prove that $\mu$ descends to a natural transformation $\mu^G\colon \Phi_1^G \to \Phi_2^G$ we 
check that for every $(c, \theta) \in \CC^G$ the morphism $\mu\colon \Phi_1(c) \to \Phi_2(c)$ commutes
with linearizations:
\[\xymatrix{
\rho_g \Phi_1(c) \ar[d]_{\delta_1}^{\simeq} \ar[r]^{\rho_g \mu(c)} & \rho_g \Phi_2(c) \ar[d]_{\delta_2}^{\simeq} \\
\Phi_1 \rho_g(c) \ar[r]^{\mu \rho_g(c)} &  \Phi_2 \rho_g(c) \\
\Phi_1 (c) \ar[u]_{\Phi_1 \theta_g} \ar[r]^{\mu(c)} &  \Phi_2(c) \ar[u]_{\Phi_2 \theta_g}\\
}\]

The transformation $\mu^G$ is natural since the original transformation $\mu$ is natural and the forgetful
functor $\CC^G \to \CC$ is faithful.

\subsection{Definition} Two weak $G$-functors $\Phi\colon \CC \to \DD$, $\Psi\colon \DD \to \CC$
are called {\bf $\mathbf{G}$-adjoint} if they are adjoint and the unit $\eps\colon id \to \Phi \Psi$
and counit $\eta\colon \Psi \Phi \to id$ of the adjunction are $G$-natural transformations.

\subsection{Lemma} 
\label{adj-descent}
A $G$-adjoint pair of functors $\Phi$, $\Psi$ induces an adjoint pair
$\Phi^G$, $\Psi^G$ between the categories of equivariant objects.

\subsection{Proof}

From \ref{lem-Gnat} it follows that we have natural transformations
$\eps^G\colon id \to \Phi^G \Psi^G$, $\eta^G\colon \Phi^G \Psi^G \to id$.
The condition for $\Psi$ and $\Phi$ to be adjoint is that two compositions
\[
\Phi \eta \circ \eps \Phi\colon \Phi \to \Phi \Psi \Phi \to \Phi 
\]
and
\[
\eta \Phi \circ \Psi \eps\colon \Psi \to \Psi \Phi \Psi \to \Psi 
\]
are identities. Since the forgetful functor $\CC^G \to \CC$ is faithful,
the same holds for $\Phi^G$, $\Psi^G$.

\subsection{Proposition}
\label{adj-mates}
A left or right adjoint $\Psi$ to a weak $G$-functor $\Phi$ can be made into a weak $G$-functor in such a way
that $\Psi$ and $\Phi$ become $G$-adjoint.
\subsection{Proof}


Let $\Psi$ be the left adjoint to $\Phi\colon \CC \to \DD$. 
We construct the structure of a left lax $G$-functor
on $\Psi$ using the structure of a right lax $G$-functor on $\Phi$.

Let $\eps\colon id \to \Phi \Psi$ and $\eta\colon \Psi \Phi \to id$ be the
unit and the counit of the adjunction.

Given a right lax $G$-structure $\delta_g\colon \rho_g \Phi \to \Phi \rho_g$ on $\Phi$
we define the left lax $G$-structure $\delta_g'\colon \Psi \rho_g \to \rho_g \Psi$ on $\Psi$ 
as a mate of $\delta_g$ with respect to the adjunction \cite[Prop. 2.1]{KS}, \cite[pp.\,185--186]{L}, i.e. 
\[
\delta_g' = \eta \rho_g \Psi \circ \Psi \delta_g \Psi \circ \Psi \rho_g \eps\colon \Psi \rho_g \to \Psi \rho_g \Phi \Psi \to \Psi \Phi \rho_g \Psi \to \rho_g \Psi.
\]

The pentagon axiom can be expressed as an equality of certain compositions in the double category of \cite[p.86]{KS}, hence
is preserved under taking mates by \cite[Prop. 2.2]{KS}. 
Checking the identity axiom for $\delta_1'$ is straightforward.

Now by \ref{lem-weak} $\Psi$ becomes a weak $G$-functor.
The proof for right adjoints is analogous.

\medskip

We now need to prove that the unit and counit transformations $\eps$, $\eta$ are $G$-natural.
We do the proof for the unit $\eps$. We need to check that the following diagram commutes:
\[\xymatrix{
\rho_g id \ar[d]_{=} \ar[r]^{\eps} & \rho_g \Phi \Psi \ar[d]^{\delta_{\Phi\Psi, g}} \\
id \, \rho_g \ar[r]^{\eps} & \Phi \Psi \rho_g. \\
}\]

Here $\delta_{\Phi\Psi}$ is defined using \ref{G-compos}. Unraveling the definitions we are left with checking
the diagram (where we use simplified notation for the natural transformations to denote the obvious compositions)
\[\xymatrix{
\rho_g \ar[r]^{\eps} \ar[d]_{\eps} & \rho_g \Phi\Psi \ar[r]^{\delta_g}_{\simeq} \ar[d]_\eps & \Phi \rho_g \Psi \ar[r]^{=} \ar[d]_{\eps} & \Phi \rho_g \Psi \\
\Phi \Psi \rho_g \ar[r]^{\eps} & \Phi\Psi\rho_g \Phi\Psi \ar[r]^{\delta_g}_{\simeq} & \Phi \Psi \Phi \rho_g \Psi \ar[ur]_\eta\\
}\]
which is easily seen to commute.

\subsection{Corollary} Let $\Phi\colon \CC \to \DD$ be a weak $G$-functor. Then the following conditions are equivalent:
\begin{enumerate}
 \item[(a)] $\Phi$ is an equivalence of categories
 \item[(b)] There exists a weak $G$-functor $\Psi\colon \DD \to \CC$ and $G$-natural isomorphisms
 $\Psi \circ \Phi \simeq \id_\CC$, $\Phi \circ \Psi \simeq \id_\DD$.
\end{enumerate}
In this case we will call $\Phi$ a {\bf weak $G$-equivalence}.

\subsection{Proof} We only need to prove $(a) \implies (b)$ as the opposite implication is trivial.
Let $\Psi\colon \DD \to \CC$ be the quasi-inverse functor to $\Phi$. In particular $\Psi$ and $\Phi$ are adjoint (both ways)
so that by \ref{adj-mates} $\Psi$ has a structure of a weak $G$-functor with compositions $G$-isomorphic to identity functors.

\section{Example: $G$-actions on the category of vector spaces}

\subsection{} In this section we review a well-known example of how equivalence classes of $G$-actions on the category of $k$-vector spaces 
correspond bijectively to cohomology classes $H^2(G, k^*)$.

\subsection{} Let $\CC = \Vect$ be the category of $k$-vector spaces, and let $\rho$ be the $G$-action on $\Vect$.
As every autoequivalence of $\CC$ is isomorphic to the identity functor, let us assume $\rho_g = \id$ for every $g \in G$.
In this setup the data of the $G$-action $\rho$ defined in \ref{def-action} is equivalent to specifying a cocycle $\phi \in Z^2(G,k^*)$.

\subsection{} Consider two $G$-actions on $\Vect$ given by cocycles $\phi, \phi' \in Z^2(G,k^*)$.
For the $G$-actions to be equivalent there needs to exist a weak $G$-functor
\[
\Phi\colon (\Vect, \phi) \to (\Vect, \phi') 
\]
which is an equivalence of categories. Then the pentagon axiom \ref{pentagon} requires existence
of an element $\delta = (\delta_g)_{g \in G} \in Z^1(G,k^*)$ such that $\phi'_{g,h} = \delta_g \delta_h \delta_{gh}^{-1} \phi_{g,h}$ for all $g,h$.
Thus $G$-categories $(\Vect, \phi)$ and $(\Vect, \phi')$ are equivalent if and only if $[\phi] = [\phi'] \in H^2(G, k^*)$.

\subsection{} The category of equivariant objects $(\Vect, \phi)^G$ is the category of $\phi$-twisted $G$-representations
with objects given by vector spaces $V$ together with isomorphism $\theta_g\colon V \to V$
satisfying $\theta_{gh} = \phi(g,h) \theta_{g} \theta_{h}$ and $G$-equivariant morphisms.
In particular, if $\phi$ is the trivial cocycle, so that $G$-action on $\Vect$ is trivial,
$\Vect^G$ is the category of $G$-representations.

\section{Strictifying $G$-actions}

\subsection{} Let $\Omega(G)$ denote the category with one object for every element $g \in G$ with $\Hom(g,g) = k$ and
$\Hom(g,h) = 0$ for $g \ne h$.

\subsection{} Let $\CC$ be a category with a $G$-action. Consider the category of weak $G$-functors and $G$-natural transformations from $\Omega(G)$ to $\CC$
\[
\CC' = \Hom_{G}(\Omega(G), \CC). 
\]
We endow $\CC'$ with the strict $G$-action induced by the $G$-action on $\Omega(G)$. 

\subsection{} Explicitly the objects of $\CC'$ consist of families $(c_g \in \CC)_{g \in G}$ together
with isomorphisms $\delta_{h,g}\colon \rho_h c_g \simeq c_{hg}$ satisfying the cocycle condition that
two ways of getting an isomorphism $\rho_k \rho_h c_g \simeq c_{khg}$ coincide.
The morphisms from $(c_g)_{g \in G}$ to $(d_g)_{g \in G}$ are morphisms $f_g\colon c_g \to d_g$
satisfying the condition that the two natural ways of forming a morphism $\rho_h c_g \to d_{hg}$ coincide.

\subsection{Theorem} \label{strict} The functor $\Phi\colon \CC' \to \CC$ sending $(c_g)_{g \in G}$ to $c_1$ is a weak $G$-equivalence.
Hence, every $G$-action is weakly equivalent to a strict $G$-action.

\subsection{Proof} We need to check that $\Phi$ has a structure of a weak $G$-functor and that $\Phi$ 
is fully faithful and essentially surjective.

The structure of a weak $G$-functor on $\Phi$ 
is in fact simply given by the structure maps $\delta_{h,g}$. That is we have
functorial isomorphisms
\[
\rho_g \Phi(c) = \rho_g(c_1) \overset{\delta_{g,1}}{\to} c_g = \Phi \rho_g(c)
\]
and the pentagon axiom follows from the cocycle condition on $\delta$.

To check that $\Phi$ is essentially surjective, one checks that for any $c \in \CC$ the
family $(\rho_g(c))$ has a structure of an object from $\CC$.
Furthermore, one can see that any object $(c_g)_{g \in G}$ is isomorphic to $(\rho_g(c_1))_{g \in G}$.

Thus to check that $\Phi$ is fully faithful, we may take two objects $(c_g)_{g \in G} = (\rho_g(c_1))$
and $(d_g)_{g \in G} = (\rho_g(d_1))$ and a morphism $f_g\colon c_g \to d_g$ between them.
It is then easy to see that $f_g = \rho(f_1)$ and that conversely for any $f_1\colon c_1 \to d_1$,
the collection $\rho_g(f_1)$ defines a morphism between $c$ and $d$.

\section{Elagin's Theorem}

\subsection{} If $\CC$ is a triangulated category and $G$ acts by triangulated autoequivalences, then $\CC^G$ is endowed
with a shift functor and a set of distinguished triangles: these are the triangles that are distinguished after applying the forgetful functor $\CC^G \to \CC$.
Furthermore under some mild technical assumptions this gives $\CC^G$ the structure of a triangulated category \cite[Theorem 6.9]{E2}, for instance
existence of a dg-enhancement of $\CC$ is a sufficient condition for $\CC^G$ to be triangulated \cite[Corollary 6.10]{E2}.

\subsection{Theorem} \label{thm-elagin}
Let $\CC = \langle \AA, \BB \rangle$ be a semi-orthogonal decomposition of triangulated categories.
Let $G$ act on $\CC$ by triangulated autoequivalences which preserve $\AA$ and $\BB$. Assume that the equivariant category $\CC^G$ is triangulated with respect
to triangles coming from $\CC$.
Then $\AA^G$, $\BB^G \subset \CC^G$ are triangulated and there is a semi-orthogonal decomposition
\[
\CC^G = \langle \AA^G, \BB^G \rangle. 
\]

\subsection{Proof}

The existence of an adjoint pair between $\CC$ and $\CC^G$ \cite[Lemma 3.7]{E2}
implies that $\BB^G = {}^\perp \AA^G$ and $\AA^G = {\BB^G}^\perp$. In particular
$\AA^G$ and $\BB^G$ are triangulated subcategories of $\CC^G$.

Now in order to establish the semi-orthogonal decomposition $\CC^G = \langle \AA^G, \BB^G \rangle$
it suffices to show that the embedding $i^G\colon \AA^G \to \CC^G$ has a left adjoint \cite[1.5]{BK1}.
This holds true by \ref{adj-mates}, \ref{adj-descent}: the functor $i\colon \AA \to \CC$ is (strictly) $G$-equivariant, hence its
left adjoint $p\colon \CC \to \AA$ induces an adjoint $p^G$ to the embedding $i^G\colon \AA^G \to \CC^G$. 



\providecommand{\arxiv}[1]{\href{http://arxiv.org/abs/#1}{\tt arXiv:#1}}

\medskip
\medskip

\medskip

\address{
{\bf Evgeny Shinder}\\
School of Mathematics and Statistics \\
University of Sheffield \\
The Hicks Building \\
Hounsfield Road \\
Sheffield S3 7RH\\
e-mail: {\tt eugene.shinder@gmail.com}
}

\end{document}